\documentclass{article}
\usepackage{amsmath,amssymb,amsthm}

\newtheorem{theorem}{Theorem}
\newtheorem{corollary}{Corollary}
\newtheorem{lemma}{Lemma}
\newtheorem{proposition}{Proposition}
\newtheorem{definition}{Definition}

\numberwithin{theorem}{section}
\numberwithin{corollary}{section}
\numberwithin{definition}{section}
\numberwithin{lemma}{section}
\numberwithin{proposition}{section}

\newcommand{\nonumsection}[1] {{\vspace{18pt}{\noindent\Large \bf #1}}\vspace{12pt}}

\begin{document}
\title{Max--Min Representation of Piecewise Linear Functions}
	\author{Sergei~Ovchinnikov \\
	Mathematics Department\\
	San Francisco State University\\
	San Francisco, CA 94132\\
	sergei@sfsu.edu} 
\date{\today} 
\maketitle

\begin{abstract}
\noindent
It is shown that any piecewise linear function can be represented as a Max--Min polynomial of its linear components.
\end{abstract}

\section{Introduction}
The goal of the paper is to establish a representation of a piecewise linear function on a closed convex domain in $\mathbb{R}^d$ as a Max--Min composition of its linear components. 

The paper is organized as follows. In Section~2 we introduce a hyperplane arrangement associated with a given piecewise linear function $f$ on a closed convex domain $\Gamma$. This arrangement defines a set $\mathcal{T}$ of regions with closures forming a cover of $\Gamma$. A distance function on $\mathcal{T}$ is introduced and its properties are established in Section~3. This distance function is an essential tool in our proof of the main result which is found in Section~4 (Theorem~\ref{main-theorem}). Some final remarks are made in Section~5.

The `standard' text on convex polytopes is~\cite{bG67}. More information on hyperplane arrangements is found in~\cite{gZ95} and~\cite{aB93}.

\section{Preliminaries}

We begin with the following definition.
\begin{definition} \label{PL-function}
Let $\Gamma$ be a closed convex domain in $\mathbb{R}^d$. A function $f:\Gamma\rightarrow\mathbb{R}$ is said to be \emph{piecewise linear} if there is a finite family $\mathcal{Q}$ of closed domains such that $\Gamma=\cup\mathcal{Q}$ and $f$ is linear on every domain in $\mathcal{Q}$. A unique linear function $g$ on $\mathbb{R}^d$ which coincides with $f$ on a given $Q\in\mathcal{Q}$ is said to be a \emph{component} of $f$.
\end{definition}

In this definition, an (affine) linear function is a function in the form 
\begin{equation*}
h(\boldsymbol{x})=\boldsymbol{a}\cdot\boldsymbol{x}+b = a_1 x_1+a_2 x_2 +\cdots + a_d x_d+b.
\end{equation*}
The equation $h(\boldsymbol{x})=0$ defines an (affine) hyperplane provided $\boldsymbol{a}\not=\boldsymbol{0}$.

Clearly, any piecewise linear function on $\Gamma$ is continuous. 

Let $f$ be a piecewise linear function on $\Gamma$ and $\{g_1,\ldots,g_n\}$ be the family of its distinct components. In what follows, we assume that $f$ has at least two distinct components.

Since components of $f$ are distinct functions, the solution set of any equation in the form $g_i(\boldsymbol{x})=g_j(\boldsymbol{x})$ for $i<j$ is either empty or a hyperplane. We denote $\mathcal{H}$ the set of hyperplanes defined by the above equations that have a nonempty intersection with the interior $\text{int}(\Gamma)$ of $\Gamma$. A simple topological argument shows that $\mathcal{H}\not=\emptyset$\,; thus $\mathcal{H}$ is an (affine hyperplane) arrangement. Let $\mathcal{T}$ be the family of nonempty intersections of the regions of $\mathcal{H}$ with $\text{int}(\Gamma)$. The elements of $\mathcal{T}$ are the connected components of $\text{int}(\Gamma)\setminus\cup\mathcal{H}$. Clearly, they are convex sets. Note that $\cup\mathcal{T}$ is dense in $\Gamma$. We shall use the same name `region' for elements of $\mathcal{T}$. The closure $\bar{Q}$ of $Q\in\mathcal{T}$ is the intersection of a polyhedron with $\Gamma$. The intersections of facets of this polyhedron with $\Gamma$ will be also called facets of $Q$. Two regions in $\mathcal{T}$ are adjacent if the intersection of their closures is a common facet of these regions.

We have the following trivial but important property of $\mathcal{T}$.

\begin{proposition}
The restrictions of the components of $f$ to any given region $Q\in\mathcal{T}$ are linearly ordered, i.e., for all $i\not=j$, either $g_i(\boldsymbol{x})>g_j(\boldsymbol{x})$ for all $\boldsymbol{x}\in Q$, or $g_i(\boldsymbol{x})<g_j(\boldsymbol{x})$ for all $\boldsymbol{x}\in Q$.
\end{proposition}

In the rest of the paper, we shall use this property of components without making explicit reference to it.

\section{Metric structure on $\mathcal{T}$}

We use a straightforward geometric approach to define a distance function on $\mathcal{T}$. It is the same distance function as in~\cite[Section~4.2]{aB93} where it is defined as the graph distance on the tope graph.

For given $P,Q\in\mathcal{T}$, let $S(P,Q)$ denote the \emph{separation set} of $P$ and $Q$, i.e., the set of all hyperplanes in $\mathcal{H}$ separating $P$ and $Q$. 

Let $\boldsymbol{p}\in P$ and $\boldsymbol{q}\in Q$ be two points in distinct regions $P$ and $Q$. Suppose $S(P,Q)=\emptyset$. Then the interval $[\boldsymbol{p},\boldsymbol{q}]$ belongs to a connected component of $\text{int}(\Gamma)\setminus\cup\mathcal{H}$ implying $P=Q$, a contradiction. Thus we may assume $S(P,Q)\not=\emptyset$.

The interval $[\boldsymbol{p},\boldsymbol{q}]$ is a subset of $\text{int}(\Gamma)$ and has a single point intersection with any hyperplane in $S(P,Q)$. Moreover, we can always choose $\boldsymbol{p}$ and $\boldsymbol{q}$ in such a way that different hyperplanes in $S(P,Q)$ intersect $[\boldsymbol{p},\boldsymbol{q}]$ in different points. Let us number these points in the direction from $\boldsymbol{p}$ to $\boldsymbol{q}$ as follows
\begin{equation*}
\boldsymbol{r}_0 = \boldsymbol{p}, \boldsymbol{r}_1,\ldots,\boldsymbol{r}_{k+1}=\boldsymbol{q}.
\end{equation*}
Each open interval $(\boldsymbol{r}_i,\boldsymbol{r}_{i+1})$ is an intersection of $[\boldsymbol{p},\boldsymbol{q}]$ with some region which we denote $R_i$ (in particular, $R_0=P$ and $R_k=Q$). Moreover, by means of this construction, points $\boldsymbol{r}_i$ and $\boldsymbol{r}_{i+1}$ belong to facets of $R_i$. We conclude that regions $R_i$ and $R_{i+1}$ are adjacent for all $i=0,1,\ldots,k-1$.

Let us define $d(P,Q)=|S(P,Q)|$ for all $P,Q\in\mathcal{T}$. It follows from the argument in the foregoing paragraph that the function $d$ satisfies the following conditions:
\begin{itemize}
	\item[(i)] $d(P,Q)=0$ if and only if $P=Q$.
	\item[(ii)] $d(P,Q)=1$ if and only if $P$ and $Q$ are adjacent regions.
	\item[(iii)] If $d(P,Q)=m$, then there exists a sequence $R_0=P,R_1,\ldots,R_m=Q$ of regions in $\mathcal{T}$ such that $d(R_i,R_{i+1})=1$ for $0\leq i<m$.
\end{itemize}

From an obvious relation $S(P,Q) = S(P,R)\,\Delta\, S(R,Q)$ it follows that 
\begin{itemize}
	\item[(iv)] $d(P,Q)\leq d(P,R)+d(R,Q)$, and
	\item[(v)] $d(P,Q)=d(P,R)+d(R,Q)$ if and only if $S(P,Q)=S(P,R)\cup S(R,Q)$.
\end{itemize}

We summarize these properties of $d$ in the following proposition.

\begin{proposition} \label{proposition}
The function $d(P,Q)=|S(P,Q)|$ is a metric on $\mathcal{T}$ satisfying the following properties:
\begin{itemize}
	\item[\emph{(i)}] $d(P,Q)=1$ if and only if $P$ and $Q$ are adjacent regions.
	\item[\emph{(ii)}] If $d(P,Q)=m$ then there exists a sequence $R_0=P,R_1,\ldots,R_m=Q$ such that $d(R_i,R_{i+1})=1$ for $0\leq i<m$.
	\item[\emph{(iii)}] $d(P,Q)=d(P,R)+d(R,Q)$ if and only if $S(P,Q)=S(P,R)\cup S(R,Q)$.
\end{itemize}
\end{proposition}

\section{Main theorem}
In this section, we prove the following theorem.

\begin{theorem} \label{main-theorem}
\emph{(a)} Let $f$ be a piecewise linear function on $\Gamma$ and $\{g_1,\ldots,g_n\}$ be the set of its distinct components. There exists a family $\{S_j\}_{j\in J}$ of subsets of $\{1,\ldots,n\}$ such that
\begin{equation} \label{polynomial}
f(\boldsymbol{x}) = \bigvee_{j\in J}\bigwedge_{i\in S_j} g_i(\boldsymbol{x}),\quad \forall \boldsymbol{x}\in\Gamma.
\end{equation}

\emph{(b)} Conversely, for any family of distinct linear functions $\{g_1,\ldots,g_n\}$ the above formula defines a piecewise linear function.
\end{theorem}
Here, $\lor$ and $\land$ are operations of maximum and minimum, respectively. The expression on the right side in~(\ref{polynomial}) is a Max--Min (lattice) polynomial in the variables $g_i$'s.

For a given $P\in\mathcal{T}$ we denote $f^P$ (resp. $g_i^P$) the restriction of $f$ (resp. $g_i$) to $P$. The functions $g_i^P$ are linearly ordered for a fixed $P$. Since the restriction of $f$ to $P$ is one of the functions $g_i^P$, there is a unique number $n(P)$ such that $f^P = g_{n(P)}^P$.

\begin{lemma} \label{lemma}
For any $P,Q\in\mathcal{T}$ there exists $k$ such that
\begin{equation*}
g_k^P \leq g_{n(P)}^P \quad\text{and}\quad g_k^Q \geq g_{n(Q)}^Q,
\end{equation*}
or, equivalently,
\begin{equation*}
g_k(\boldsymbol{x}) \leq f(\boldsymbol{x}),\;\forall\boldsymbol{x}\in P,\quad\text{and}\quad g_k(\boldsymbol{x}) \geq f(\boldsymbol{x}),\;\forall\boldsymbol{x}\in Q.
\end{equation*}
\end{lemma}

\begin{proof}
The proof is by induction on $d(P,Q)$.

(i) $d(P,Q)=1$. By Proposition~\ref{proposition}\,(i), $P$ and $Q$ are adjacent regions. Let $F$ be the common facet of the closures of $P$ and $Q$ and let $H$ be the affine span of $F$. Since functions $f,g_{n(P)},g_{n(Q)}$ are continuous, $g_{n(P)}(\boldsymbol{x})=g_{n(Q)}(\boldsymbol{x})$ for all $\boldsymbol{x}\in F$ and therefore for all $\boldsymbol{x}\in H$. We may assume that $g_{n(P)}^P<g_{n(Q)}^P$ (the other case is treated similarly). Then $g_{n(P)}^Q>g_{n(Q)}^Q$, and $k=n(P)$ satisfies conditions of the lemma.

(ii) $d(P,Q)>1$. By Proposition~\ref{proposition}\,(ii) and (i), there is a region $R$ adjacent to $P$ such that $d(R,Q)=d(P,Q)-1$. By the induction hypothesis, there is $r$ such that
\begin{equation*}
g_r^R \leq g_{n(R)}^R \quad\text{and}\quad g_r^Q \geq g_{n(Q)}^Q.
\end{equation*}
If $g_r^P\leq g_{n(P)}^P$, then $k=r$ satisfies conditions of the lemma. Otherwise, we have $g_r^P > g_{n(P)}^P$. By Proposition~\ref{proposition}\,(iii), the unique hyperplane $H\in\mathcal{H}$ that separates $P$ and $R$ also separates $P$ and $Q$. The same argument as in (i) shows that $g_{n(P)}(\boldsymbol{x})=g_{n(R)}(\boldsymbol{x})$ for all $\boldsymbol{x}\in H$. Since $g_r^P > g_{n(P)}^P$ and $g_r^R \leq g_{n(R)}^R$, we have $g_r(\boldsymbol{x})=g_{n(P)}(\boldsymbol{x})$ for all $\boldsymbol{x}\in H$. Consider function $g=g_r-g_{n(P)}$. It is zero on the hyperplane $H$ and positive on the full--dimensional region $P$. Thus it is positive on the open halfspace containing $P$. Hence, it must be negative on the open halfspace containing $R$ and $Q$. We conclude that $g_{n(P)}^Q > g_r^Q$. Since $g_r^Q \geq g_{n(Q)}^Q$, $k=n(P)$ satisfies conditions of the lemma.

\end{proof}

Now we proceed with the proof of Theorem~\ref{main-theorem}.

(a) For a given $P\in\mathcal{T}$ we define $S_P=\{i:g_i^P\geq g_{n(P)}^P\}\subseteq\{1,\ldots,n\}$. Let $F_P(\boldsymbol{x})$ be a function defined on $\Gamma$ by the equation
\begin{equation*}
F_P(\boldsymbol{x}) = \bigwedge_{i\in S_P} g_i(\boldsymbol{x}).
\end{equation*}
Clearly, $F_P(\boldsymbol{x})=g_{n(P)}(\boldsymbol{x})=f(\boldsymbol{x})$ for all $\boldsymbol{x}\in P$.

Suppose $F_P(\boldsymbol{y}) < F_Q(\boldsymbol{y})$ for some $\boldsymbol{y}\in P$ and $Q\not=P$, i.e.,
\begin{equation*}
\bigwedge_{i\in S_P} g_i(\boldsymbol{y}) < \bigwedge_{j\in S_Q} g_j(\boldsymbol{y}).
\end{equation*}
Then $g_{n(P)}^P < g_j^P$ for all $j\in S_Q$. Thus for any $j$ such that $g_j^Q\geq g_{n(Q)}^Q$ we have $g_{n(P)}^P < g_j^P$. This contradicts Lemma~\ref{lemma}. 

Hence, $F_Q(\boldsymbol{x})\leq F_P(\boldsymbol{x})=f(\boldsymbol{x})$ for all $\boldsymbol{x}\in P$ and $Q\in\mathcal{T}$.

Consider function $F(\boldsymbol{x})$ defined by
\begin{equation*}
F(\boldsymbol{x}) = \bigvee_{P\in\mathcal{T}}F_P(\boldsymbol{x}) = \bigvee_{P\in\mathcal{T}}\bigwedge_{i\in S_P} g_i(\boldsymbol{x})
\end{equation*}
for all $\boldsymbol{x}\in\Gamma$. Clearly, $f(\boldsymbol{x})=F(\boldsymbol{x})$ for all $\boldsymbol{x}\in\cup\,\mathcal{T}$. Since $\cup\,\mathcal{T}$ is dense in $\Gamma$ and $f$ and $F$ are continuous functions on $\Gamma$, we conclude that
\begin{equation*}
f(\boldsymbol{x}) = \bigvee_{P\in\mathcal{T}}\bigwedge_{i\in S_P} g_i(\boldsymbol{x})
\end{equation*}
for all $\boldsymbol{x}\in\Gamma$. 

(b) Let $\{g_1,\ldots,g_n\}$ be a family of distinct linear functions on $\Gamma$ and let $f$ be defined by~(\ref{polynomial}). Consider sets $H_{ij}=\{\boldsymbol{x}:g_i(\boldsymbol{x})=g_j(\boldsymbol{x}),\;i>j\}$. If the intersections of these sets with $\text{int}(\Gamma)$ are empty, then the functions $g_i$'s are linearly ordered over $\text{int}(\Gamma)$ and, by~(\ref{polynomial}), $f$ is a linear function. Otherwise, let $\mathcal{H}$ be the family of sets $H_{ij}$'s with nonempty intersections with $\text{int}(\Gamma)$. Let $Q$ be a region of the arrangement $\mathcal{H}$. Since $Q$ is connected, the functions $g_i$'s are linearly ordered over $Q$ and, by~(\ref{polynomial}), there is $i$ such that $f(\boldsymbol{x})=g_i(\boldsymbol{x})$ for all $\boldsymbol{x}\in Q$. The same is also true for the closure of $Q$.

This completes the proof of Theorem~\ref{main-theorem}.

\begin{corollary}
Let $\Gamma$ be a star--like domain in $\mathbb{R}^d$ such that its boundary $\partial\Gamma$ is a polyhedral complex. Let $f$ be a function on $\partial\Gamma$ such that its restriction to each $(d-1)$-dimensional polyhedron in $\partial\Gamma$ is a linear function on it. Then $f$ admits representation~\emph{(\ref{polynomial})}.
\end{corollary}

\begin{proof}
Let $\boldsymbol{a}$ be a central point in $\Gamma$. For $\boldsymbol{x}\in\mathbb{R}^d,\;\boldsymbol{x}\not=\boldsymbol{a}$, let $\tilde{\boldsymbol{x}}$ be the unique intersection point of the ray from $\boldsymbol{a}$ through $\boldsymbol{x}$ with $\partial\Gamma$. We define
\begin{equation*}
\tilde{f}(\boldsymbol{x}) = \begin{cases}
	\frac{\Vert\boldsymbol{x}-\boldsymbol{a}\Vert}{\Vert\tilde{\boldsymbol{x}}-\boldsymbol{a}\Vert}f(\tilde{\boldsymbol{x}}), &\text{for $\boldsymbol{x}\not=\boldsymbol{a}$,} \\
	0, &\text{for $\boldsymbol{x}=\boldsymbol{a}$.}
\end{cases}
\end{equation*}
Clearly, $\tilde{f}$ is a piecewise linear function on $\mathbb{R}^d$ and $\tilde{f}\vert_{\partial\Gamma} = f$. Thus $f$ admits representation~(\ref{polynomial}).

\end{proof}

Note that the previous corollary holds for any polyhedron in $\mathbb{R}^d$.

\section{Concluding remarks}
\begin{enumerate}
	\item The statements of Theorem~\ref{main-theorem} also hold for piecewise linear functions from $\Gamma$ to $\mathbb{R}^m$. Namely, let $\boldsymbol{f}:\Gamma\rightarrow\mathbb{R}^m$ be a piecewise linear function and $\{\boldsymbol{g}_1,\ldots,\boldsymbol{g}_n\}$ be the set of its distinct components. We denote 
\begin{equation*}
\boldsymbol{f}=(f_1,\ldots,f_m)\quad\text{and}\quad\boldsymbol{g}_k=(g_1^{(k)},\ldots,g_m^{(k)})\text{,~for~}1\leq k\leq n.
\end{equation*}
There exists a family $\{S_j^k\}_{j\in J,\;1\leq k\leq n}$ of subsets of $\{1,\ldots,n\}$ such that
\begin{equation*}
f_k(\boldsymbol{x}) = \bigvee_{j\in J}\bigwedge_{i\in S_j^k} g_i^{(k)}(\boldsymbol{x}),\quad \forall \boldsymbol{x}\in\Gamma,\;1\leq k\leq m.
\end{equation*}
The converse is also true.
	\item The convexity of $\Gamma$ is an essential assumption. Consider, for instance, the domain in $\mathbb{R}^2$ which is a union of three triangles defined by the sets of their vertices as follows:
\begin{align*}
&\Delta_1=\{(-1,0),(-1,-1),(0,0)\},\quad\Delta_2=\{(0,0),(1,1),(1,0)\}, \\
& \text{and}\quad\Delta_3=\{(-1,0),(1,0),(0,-1)\}.
\end{align*}
Let us define
\begin{equation*}
f(\boldsymbol{x}) = \begin{cases}
	-x_2, &\text{for $\boldsymbol{x}\in\Delta_1,$} \\
	x_2, &\text{for $\boldsymbol{x}\in\Delta_2\cup\Delta_3,$}
\end{cases}
\end{equation*}
where $\boldsymbol{x}=(x_1,x_2)$. This piecewise linear function has two components, $g_1(\boldsymbol{x})=-x_2$ and $g_2(\boldsymbol{x})=x_2$, but is not representable in the form~(\ref{polynomial}).
	\item Likewise,~(\ref{polynomial}) is not true for piecewise polynomial functions as the following example (due to B.~Sturmfels) illustrates. Let $\Gamma=\mathbb{R}^1$. We define
\begin{equation*}
f(x) = \begin{cases}
	0, &\text{for $x\leq 0$,} \\
	x^2, &\text{for $x>0$.}
\end{cases}
\end{equation*}
	\item It follows from Theorem~\ref{main-theorem} that any piecewise linear function on a closed convex domain in $\mathbb{R}^d$ can be extended to a piecewise linear function on the entire space $\mathbb{R}^d$.
\end{enumerate}

\nonumsection{Acknowledgments}

\noindent
The author thanks S.~Gelfand, O.~Musin, and B.~Sturmfels for helpful discussions on the earlier versions of the paper.


\begin{thebibliography}{9}
\bibitem{aB93}
	A.~Bj\"{o}rner, M.~Las~Vergnas, B.~Sturmfels, N.~White, and G.M.~Ziegler, Oriented Matroids, Encyclopedia of Mathematics, Vol.~46, Cambridge University press, 1993.
\bibitem{bG67}
	B.~Gr\"{u}nbaum, Convex Polytopes, Interscience, London, 1967.
\bibitem{gZ95}
	G.M.~Ziegler, Lectures on Polytopes, Graduate Text in Mathematics, Springer--Verlag, 1995.
\end{thebibliography}
\end{document}